\documentclass[12pt]{amsart}
\usepackage{amsfonts,amsthm,latexsym,amsmath,amssymb,amscd,amsmath, mathrsfs, url, cite }
\usepackage[pdftex]{color}
\usepackage[bookmarks=true,hyperindex,pdftex,colorlinks, citecolor=blue,linkcolor=blue, urlcolor=blue]{hyperref}
\newcounter{count}
\numberwithin{count}{section}
\newtheorem{Lemma}[count]{Lemma}

\newtheorem{Example}[count]{Example}
\newtheorem{Definition}[count]{Definition}

\newtheorem{Theorem}[count]{Theorem}
\newtheorem{Conjecture}[count]{Conjecture}
\newtheorem{Problem}[count]{Problem}
\newtheorem{Statement}[count]{Statement}

\begin{document}

\author[O.~Katkova]{Olga Katkova}

\address{Department of Mathematics, University of Massachusetts Boston, USA}
\email{olga.katkova@umb.edu }

\author[A.~Vishnyakova]{Anna Vishnyakova}

\address{Department of Mathematics, Holon Institute of Technology,
Israel}
\email{annalyticity@gmail.com}

\title[Multiplier sequences for the totally positive sequences]
{An analog of multiplier sequences for the set of totally positive sequences}

\begin{abstract}
A real sequence
$(b_k)_{k=0}^\infty$  is called totally positive  
if all minors of the infinite matrix $ \left\| b_{j-i}  
\right\|_{i, j =0}^\infty$ are nonnegative (here $b_k=0$
for $k<0$). In this paper, we  investigate the problem
of description of the set of sequences $(a_k)_{k=0}^\infty$
such that for every totally positive sequence $(b_k)_{k=0}^\infty$ 
the sequence $(a_k b_k)_{k=0}^\infty$ is also totally positive.
We obtain the description of such sequences $(a_k)_{k=0}^\infty$
in two cases: when the generating function of the sequence
$\sum_{k=0}^\infty a_k z^k$ has at least one pole, and when
the sequence $(a_k)_{k=0}^\infty$ has not more than $4$
nonzero terms.

\end{abstract}



\maketitle

\section{Introduction}

\medskip
\begin{center}
\textit{To our dear teacher Iossif Vladimirivich Ostrovskii with love and gratitude. 
We will always remember and admire him not only as
an outstanding mathematician and a brilliant teacher, but also as a person 
of amazing kindness, generosity and nobility.}
\end{center}
\medskip

We start with the definition of multiply positive and totally positive sequences.

\begin{Definition}
A real sequence $(a_k)_{k=0}^\infty$  is called  $m$-times positive 
($m \in \mathbb{N}$), if all minors of the
infinite matrix
\begin {equation}
\label{mat}
 \left\|
  \begin{array}{ccccc}
   a_0 & a_1 & a_2 & a_3 &\ldots \\
   0   & a_0 & a_1 & a_2 &\ldots \\
   0   &  0  & a_0 & a_1 &\ldots \\
   0   &  0  &  0  & a_0 &\ldots \\
   \vdots&\vdots&\vdots&\vdots&\ddots
  \end{array}
 \right\|
\end {equation}
of orders less than or equal to $m$ are nonnegative. The class of
$m$-times positive sequences is denoted by $TP_m,$ the class of 
the generating functions of
$m$-times positive sequences ($f(x) = \sum_{k=0}^\infty a_k x^k$)
is denoted by $\widetilde{TP}_m.$
\end{Definition}

\begin{Definition}
A real sequence
$(a_k)_{k=0}^\infty$  is called totally positive  
if all minors of the infinite matrix
(\ref{mat}) are nonnegative.
The class of totally positive sequences is denoted by $TP_\infty,$
the class of the generating functions of
totally positive sequences is denoted by $\widetilde{TP}_\infty.$
\end{Definition}

Multiply positive sequences (also called P\'olya frequency
sequences) were introduced by Fekete in 1912 (see~\cite{fek}) in
connection with the problem of exact calculation of the number of
positive zeros of a real polynomial. Multiply positive and totally 
positive sequences arise in many areas of mathematics and its 
applications, see, for example, \cite{ando}, \cite{tp} or \cite{pin}.           

The class $TP_\infty$ was completely described by Aissen,
Schoenberg, Whitney and Edrei in~\cite{aissen} (see also
\cite[p. 412]{tp}).

{\bf Theorem ASWE. } {\it A function $f \in \widetilde{TP}_\infty$ if 
and only if
\begin{equation}
\label{aswe}
f(z)=C z^q e^{\gamma z}\prod_{k=1}^\infty \frac{ (1+\alpha_kz)}{(1-\beta_kz)},
\end{equation}
where $C\ge 0, q \in \mathbb{Z},
\gamma\ge 0,\alpha_k\ge 0,\beta_k\ge 0, \sum_{k=1}^\infty(\alpha_k+\beta_k)
<\infty.$ }

Theorem ASWE gives the description of the class  $TP_\infty$ in terms of 
independent  parameters  $C, q, \gamma,\alpha_k,\beta_k.$  It is easy to see that 
the class $TP_2$ consists of the sequences $(a_k)_{k=0}^\infty$ 
of the form $a_n = e^{- \psi(n)},$  where $\psi : \mathbb{N}\cup \{0 \} \to (-\infty,
 + \infty]$ is a convex function.   In \cite{ivzh} the description of the subclass of $TP_3,$  
 which consists of the sequences whose all sections belong to $TP_3,$ in terms of 
 independent parameters was obtained. The problem of the description of the
classes $TP_m$, $3 \leq m < \infty,$ in terms of independent parameters has not
been solved until now.

By theorem ASWE a polynomial $p(z) = \sum _{k=0}^n a_k z^k, \  a_k
\geq 0 ,$ has only real zeros if and only if  $(a_0,
a_1, \ldots , a_n, 0, 0, \ldots ) \in TP_\infty .$

In general, the problem of understanding whether a given polynomial has only 
real zeros is not trivial, often such problems  are  very difficult.  However, in 1926, 
J. I. Hutchinson found the following simple 
sufficient condition in terms of coefficients for an entire function with positive coefficients to 
have only real zeros.

{\bf Theorem A} (J. ~I. ~Hutchinson, \cite{hut}). { \it Let $f(x)=
\sum_{k=0}^\infty a_k x^k$, $a_k > 0$ for all $k$. 
Then $\frac{a_{n-1}^2}{a_{n-2}a_n}\geq 4$ for all $n\geq 2,$  
if and only if the following two conditions are fulfilled:\\
(i) The zeros of $f(x)$ are all real, simple and negative, and \\
(ii) the zeros of any polynomial $\sum_{k=m}^n a_k x^k$, $m < n,$  formed 
by taking any number 
of consecutive terms of $f(x) $, are all real and non-positive.}

For some extensions of Hutchinson's results see,
for example, \cite[\S4]{cc1} and \cite{HishAn}.

The question about  whether or not a given polynomial has only real zeros is 
of great importance in many areas of mathematics. So, the problem to describe 
the set of operators that preserve this set of polynomials is of the great interest.  
In connection with this problem, we define  multiplier 
sequences.

\begin{Definition} A sequence $(\gamma_k)_{k=0}^\infty$ of real
numbers is called a multiplier sequence if, whenever a real
polynomial $P(x) = \sum_{k=0}^n a_k z^k $ has only real zeros, the
polynomial $\sum_{k=0}^n \gamma_k a_k z^k $ has only real zeros. 
The class of multiplier sequences is denoted by $\mathcal{MS}$.
\end{Definition}

A simple example of a multiplier sequence is the following
sequence: $\gamma_k =k, k=0, 1, 2, \ldots$ For an arbitrary 
polynomial $P(x) = \sum_{k=0}^n a_k z^k $ with real coefficients
and only real zeros we have $\sum_{k=0}^n k a_k z^k = z P^\prime(z),$ 
and this polynomial obviously also has only real zeros.

The full description of the set of multiplier sequences was given by G.~P\'olya and 
J.~Schur in 1914. To formulate this famous result, we need a notion of
 the Laguerre-P\'olya class of  entire functions.

\begin{Definition} 
A real entire function $f$ is said to be in the {\it Laguerre--
P\'olya class of type I}, 
written $f \in \mathcal{L-P} I$, if it can
be expressed in the following form  
\begin{equation}  \label{lpc1}
 f(z) = c z^n e^{\beta z}\prod_{k=1}^\infty
\left(1+\frac {z}{x_k} \right),
\end{equation}
where $c \in  \mathbb{R},  \beta \geq 0, x_k >0 $, 
$n$ is a nonnegative integer,  and $\sum_{k=1}^\infty x_k^{-1} <
\infty$.  
\end{Definition}

Note that the product on the right-hand side can be
finite or empty (in the latter case, the product equals 1).

This class is essential in the theory of entire functions since it appears that the polynomials with 
only real and nonpositive zeros converge locally 
uniformly to these and only these functions. The following  prominent theorem
 provides an even stronger result. 

{\bf Theorem B} (E.~Laguerre and G.~P\'{o}lya, see, for example,
\cite[p. ~42--46]{HW}) and \cite[chapter VIII, \S 3]{lev}). {\it   

(i) Let $(P_n)_{n=1}^{\infty},\  P_n(0)=1, $ be a sequence
of real polynomials having only real negative zeros which  
converges uniformly on the disc $|z| \leq A, A > 0.$ Then this 
sequence converges locally uniformly in $\mathbb{C}$ to an entire function
 from the class $\mathcal{L-P}I.$
 
(ii) For any $f \in \mathcal{L-P}I$ there is a
sequence of real polynomials with only real nonpositive 
zeros, which converges locally uniformly to $f$.}

The following theorem fully describes multiplier sequences.

{\bf Theorem C} (G. P\'olya and J.Schur, cf. \cite{polsch}, \cite[pp. 100-124]{pol}
and\cite[pp. 29-47]{O}).  
{\it Let  $(\gamma_k)_{k=0}^\infty$ be a given 
real sequence. The following 
three statements are equivalent.}

{\it

1. $(\gamma_k)_{k=0}^\infty$ is a multiplier sequence.

2.  For every $n\in \mathbb{N}$ the
polynomial $P_n(z) =\sum_{k=0}^n {\binom{n}{k}} \gamma_k z^k $ has only real 
zeros of the same sign. 

3. The power series $ \Phi (z) := \sum_{k=0}^\infty \frac
{\gamma_k}{k!}z^k$ converges absolutely in the whole complex plane
and the entire function $\Phi(z)$ or the entire function
$\Phi(-z)$ admits the representation
\begin{equation}
\label{pred}
 c  z^n e^{\beta z} \prod_{k=1}^\infty
\left(1+\frac{z}{x_k}\right),
\end{equation}

where $c\in{\mathbb{R}}, \beta \geq 0, n\in {\mathbb{N}}\cup \{0\}, 0<x_k
\leq \infty,\   \sum_{k=1}^\infty \frac{1}{x_k} < \infty.$}

Strikingly, the following fact is an obvious consequence.

{\bf Corollary of Theorem C.}   {\it The sequence 
$(\gamma_0, \gamma_1, \ldots, \gamma_l,$ $ 0, 0, \ldots )$ is a
multiplier sequence if and only if the polynomial $P(z)=
\sum_{k=0}^l \frac {\gamma_k}{k!}z^k$ has only real zeros of the
same sign.}

As we mentioned before, the set of  polynomials with nonnegative
coefficients having only real nonpositive roots is a subset of the set 
$\widetilde{TP}_\infty.$  In this paper, we discover an analog of the
multiplier sequences for the set of totally positive sequences. To
formulate the problem, we need the next definition.

\begin{Definition}  Let ${\mathbf A} =(a_k)_{k=0}^\infty$ be  a nonnegative sequence.
We define the following linear convolution operator on the set of
real sequences:
$$\Lambda_{\mathbf A} ( (b_k)_{k=0}^\infty) =  (a_k b_k)_{k=0}^\infty. $$
\end{Definition} 

The following problem was posed by Alan Sokal during the 
inspiring AIM workshop ``Theory and applications of total positivity'',
July 24-July 28, 2023 (see \cite{conf} for more details). 

\begin{Problem} \label{P1}  To describe the set of  nonnegative sequences 
${\mathbf A} =(a_k)_{k=0}^\infty,$ such that the corresponding convolution 
operator $\Lambda_{\mathbf A}$ preserves the set of $TP_\infty-$sequences:
for every  $(b_k)_{k=0}^\infty \in TP_\infty$ we have $\Lambda_{\mathbf A} ( 
(b_k)_{k=0}^\infty) \in  TP_\infty. $  
\end{Problem}

For some questions connected with the problem above see \cite{DyaSok} by
A.~Dyachenko and A.~Sokal (see also previous works of A.~Dyachenko \cite{Dya1} 
and \cite{Dya2}).

We consider the multiplier sequence ${\mathbf \Gamma} = (k)_{k=0}^\infty$ and 
the corresponding convolution operator $\Lambda_{\mathbf \Gamma} 
( (b_k)_{k=0}^\infty) =  (k b_k)_{k=0}^\infty. $ As we mentioned earlier, 
this operator preserves the set of finite totally positive sequences (in other words, 
the set of coefficients of polynomials with  nonnegative coefficients and only real zeros).
But this operator does not preserve the set of all totally positive sequences.
Indeed, let us consider the function $f(z) = \frac{1}{(1-z)(2-z)}= \sum_{k=0}^\infty
b_k z^k$  (we have  $b_k =1 - \frac{1}{2^{k+1}}$). By Theorem ASWE, $ (b_k)_{k=0}^\infty
\in TP_\infty.$ But $\sum_{k=0}^\infty k b_k z^k = z f^\prime(z) =
 \frac{z(3 - 2z)}{(1-z)^2 (2-z)^2},$ this function has a positive zero, so the sequence
 of its coefficients is not a $TP_\infty-$sequence.

We will denote by $A$ the generating function of a sequence 
 ${\mathbf A} =(a_k)_{k=0}^\infty:$ $A(z) =\sum_{k=0}^\infty a_k z^k.$

Suppose that the sequence ${\mathbf A}$ has the property that the corresponding convolution 
operator $\Lambda_{\mathbf A}$ preserves the set of $TP_\infty-$sequences. Then,
since the constant sequence of all ones is the $TP_\infty-$sequence, by theorem ASWE,
the generating function $A(x)$ is a meromorphic function having the representation (\ref{aswe}). 
The following theorem gives the full description of the generating functions of 
$TP_\infty$-preservers that have at least one pole.

\begin{Theorem} \label{tt1}  Let ${\mathbf A} =(a_k)_{k=0}^\infty$ 
be a nonnegative sequence, and its generating function is a meromorphic function with
at least one pole. Then  for every  $(b_k)_{k=0}^\infty \in TP_\infty$ we have  $\Lambda_{\mathbf A} 
((b_k)_{k=0}^\infty) \in  TP_\infty$  if and only if $A(z) = \frac{C}{1- \beta z}, C > 0, \beta > 0.$  
\end{Theorem}

It remains to describe $TP_\infty$-preservers which generating functions are entire 
functions. We start with an obvious case of one or two term sequences. Let us consider 
a nonnegative sequence  ${\mathbf A} =(a_k)_{k=0}^\infty,$ 
such that $a_0 \geq 0, a_1 \geq 0,$   and $a_k=0$ for $k\geq 2.$ Then, obviously, 
for every  $(b_k)_{k=0}^\infty \in TP_\infty$  
we have  $\Lambda_{\mathbf A} ((b_k)_{k=0}^\infty) \in  TP_\infty.$ The case of
three term sequences is also simple. The following statement is obvious.

\begin{Statement} \label{St1}  Let ${\mathbf A} =(a_k)_{k=0}^\infty$ 
be a nonnegative sequence, such that $a_k>0$ for $k = 0, 1, 2, $   and 
$a_k=0$ for $k\geq 3.$ Then  for every  $(b_k)_{k=0}^\infty \in TP_\infty$ we have  
$\Lambda_{\mathbf A} ((b_k)_{k=0}^\infty) \in  TP_\infty$  if and only if 
$A(z) = a_0 + a_1 z + a_2 z^2 $  has only real (and negative) zeros. 
Moreover, $\Lambda_{\mathbf A}:  TP_\infty \to TP_\infty $ if and
only if  $\Lambda_{\mathbf A}:  TP_2 \to TP_\infty. $
\end{Statement}

The following theorem gives the description of $TP_\infty$-preservers which 
generating functions are polynomials of degree $3.$

\begin{Theorem} \label{tt2}  
 Let ${\mathbf A} =(a_k)_{k=0}^\infty$ 
be a nonnegative sequence, such that $a_k>0$ for $0\leq k \leq 3,$  and 
$a_k=0$ for $k\geq 4.$ Then  for every  $(b_k)_{k=0}^\infty \in TP_\infty$  
we have  $\Lambda_{\mathbf A} ((b_k)_{k=0}^\infty) \in  TP_\infty$  if and only 
if both polynomials $\sum_{k=0}^3 a_k x^k$ and $\sum_{k=1}^3 a_k x^k$
have only real (and nonpositive) zeros. Moreover, $\Lambda_{\mathbf A}:  TP_\infty \to TP_\infty $ if and
only if  $\Lambda_{\mathbf A}:  TP_3 \to TP_\infty. $
\end{Theorem}

Using the methods analogous to those that was used in the proof of 
Theorem~\ref{tt2}, we can prove the following statement.

\begin{Theorem} \label{tt3}  
 Let ${\mathbf A} =(a_k)_{k=0}^\infty$ 
be a nonnegative sequence, such that $a_k>0$ for $0\leq k \leq 4,$  and 
$a_k=0$ for $k\geq 5.$ Then  for every  $(b_k)_{k=0}^\infty \in TP_\infty$  
we have  $\Lambda_{\mathbf A} ((b_k)_{k=0}^\infty) \in  TP_\infty$  if and only 
if three polynomials $\sum_{k=0}^4 a_k x^k,$  $\sum_{k=1}^4 a_k x^k$
and   $\sum_{k=2}^4 a_k x^k$ have only real (and nonpositive) zeros. 
Moreover, $\Lambda_{\mathbf A}:  TP_\infty \to TP_\infty $ if and
only if  $\Lambda_{\mathbf A}:  TP_4 \to TP_\infty. $
\end{Theorem}

We will not present the proof of the above result here, since it is very cumbersome 
and does not provide a complete solution to the problem of the description of all 
entire $TP_\infty$-preservers.

The following example was given by Alan Sokal.

\begin{Example} \label{E1}  
Let $f$ be an entire function of the form $f(z) = \sum_{k=0}^\infty
a_k z^k$ with $a_0 = a_1 =1,$ $a_k =\frac{1}{q_2^{k-1}q_3^{k-2}
\cdot \ldots \cdot q_{k-1}^2 q_{k}}$ for $k\geq 2,$ where
$(q_k)_{k=2}^\infty$ is a sequence of arbitrary parameters 
under the following conditions: $q_k \geq 4$ for all $k.$
Suppose that  $(b_k)_{k=0}^\infty \in TP_\infty$   is an arbitrary 
sequence. For an entire function $ A\ast B (z) =\sum_{k=0}^\infty 
a_k b_k z^k $ we have $\frac{(a_{n-1} b_{n-1})^2}{(a_{n-2} b_{n-2)}
(a_n b_n)} = \frac{a_{n-1}^2}{a_{n-2}a_n} \cdot 
\frac{b_{n-1}^2}{b_{n-2}b_n}\geq 4$ for all $n\geq 2,$  
since $\frac{a_{n-1}^2}{a_{n-2}a_n}= q_n \geq 4$ by our assumption, 
and $\frac{b_{n-1}^2}{b_{n-2}b_n}\geq 1,$ because every 
$TP_\infty-$sequence is, in particular, $2-$times positive sequence.
Thus, using Theorem~A by Hutchinson, we get
$ A\ast B (z) \in TP_\infty.   $
\end{Example}

We formulate the following conjecture, which is consistent with 
Theorem~\ref{tt2}, Theorem~\ref{tt3} and Example~\ref{E1}.

\begin{Conjecture} \label{C1}  
 Let ${\mathbf A} =(a_k)_{k=0}^\infty$ be a nonnegative sequence. 
 Then this sequence is a $TP_\infty$-preserver, i.e.  for every  
 $(b_k)_{k=0}^\infty \in TP_\infty$  we have  
 $\Lambda_{\mathbf A} ((b_k)_{k=0}^\infty) \in  TP_\infty$  if and only 
if for every $l \in \mathbb{N} \cup \{0\}$ a formal power series $\sum_{k=l}^\infty
a_k z^k$ is an entire function from the $ \mathcal{L-P} I$ class (in particular, it has only 
real nonpositive zeros). 
\end{Conjecture}

We note that entire functions  whose Taylor sections have only real zeros 
were studied in various works (see, for example, \cite{kosshap} or \cite{klv1}), 
but entire functions whose remainders have only real zeros have been studied 
less (some results could be found in the very interesting survey \cite{iv}).
We mention here a way to construct such a function. 
The entire function $g_a(z) =\sum _{j=0}^{\infty} z^j a^{-j^2}$, $a>1,$
is called the \textit{partial theta-function}. The  survey \cite{War} by S.~O.~Warnaar 
contains the history of investigation of the partial theta-function and some of its 
main properties.  The paper \cite{klv} answers to the question:  for which 
$a>1$ the function $g_a$ belongs to the class $\mathcal{L-P} I. $ In particular,
in \cite{klv} it is proved that there exists a constant $q_\infty \approx 3{.}23363666 \ldots,  $
such that $g_a \in \mathcal{L-P} I $ if and only if $a^2 \geq q_\infty.$ In \cite{ngthv1}
the following theorem is proved. Let  $f(z) = \sum_{k=0}^\infty
a_k z^k$ with $a_0 = a_1 =1,$ $a_k =\frac{1}{q_2^{k-1}q_3^{k-2}
\cdot \ldots \cdot q_{k-1}^2 q_{k}}$ for $k\geq 2,$ where $(q_k)_{k=2}^\infty$ 
is a sequence of arbitrary parameters under the following conditions:
$q_2 \geq q_3 \geq q_4 \geq \ldots$ and $\lim_{n\to\infty}q_n \geq q_\infty.$
Then $f \in \mathcal{L-P} I.$ Using this theorem we conclude that such entire
function $f$ has all remainders with only real zeros.

\section{Proof of Theorem \ref{tt1}}

Suppose at first that $A(z) = \frac{C}{1- \beta z}, C > 0, \beta > 0.$  
Then we have $A(z) =\sum_{k=0}^\infty C\beta^k z^k,$ whence
for every  $ {\bf B} =(b_k)_{k=0}^\infty$   $ \in TP_\infty$ with  the 
generation function $B,$ the generation function of $ \Lambda_{\mathbf A} 
((b_k)_{k=0}^\infty) $ is equal to $ \sum_{k=0}^\infty C\beta^k b_k z^k =
C B(\beta  z) \in \widetilde{TP_\infty}. $ The sufficiency is proved.

Let us prove necessity. Let ${\mathbf A} =(a_k)_{k=0}^\infty$ be  such a sequence that the 
corresponding convolution operator $\Lambda_{\mathbf A}$ preserves the set 
of $TP_\infty$-sequences, and ${\mathbf A}$ is not identical zero.

\begin{Definition} 
For a nonnegative sequence ${\mathbf A} =(a_k)_{k=0}^\infty$ with 
the generating function  $A(z) =\sum_{k=0}^\infty a_k z^k$ and
a nonnegative sequence   ${\mathbf B} =(b_k)_{k=0}^\infty$ with the
generating function $B(z) =\sum_{k=0}^\infty b_k z^k$ we will denote
by $A\ast B$ the generating function of the sequence $\Lambda_{\mathbf A} 
( (b_k)_{k=0}^\infty):$
$$ (A\ast B)(z) = \sum_{k=0}^\infty a_k b_k z^k. $$ 
\end{Definition}

We mention some simple properties of the generating functions of $TP_\infty$-preservers.

\begin{Lemma} \label{L1} Suppose that a sequence ${\mathbf A} =(a_k)_{k=0}^\infty$ 
is such that  for every  $(b_k)_{k=0}^\infty \in TP_\infty$ we have $\Lambda_{\mathbf A} 
( (b_k)_{k=0}^\infty) \in  TP_\infty. $  Then
the following is true
\begin{enumerate}
\item $A(z) \in \widetilde{TP_\infty}.$

\item $A^\prime(z) \in \widetilde{TP_\infty}.$

\item $(zA(z))^\prime \in \widetilde{TP_\infty}.$

\item $\frac{1}{1 - c}(A(z) - c A(c z)) \in \widetilde{TP_\infty}$
for all $c \in (0, 1) \cup (1, \infty).$

\item $\frac{1}{1 - c}(A(z) - A(c z)) \in \widetilde{TP_\infty}$
for all $c \in (0, 1) \cup (1, \infty).$

\item $(A(z)- \frac{d}{d +1}a_0) \in \widetilde{TP_\infty}$
for all $d \geq 0.$

\item For all $n \in \mathbb{N} \cup \{0\}$ we have
 $(A(z)- \sum_{k=0}^n a_k z^k) \in \widetilde{TP_\infty}.$

\end{enumerate}
\end{Lemma}
{\it Proof of Lemma \ref{L1}.} 
\begin{enumerate}

\item We choose the sequence ${\mathbf B} =(b_k)_{k=0}^\infty  \in TP_\infty$ 
such that $b_k \equiv 1, B(z) = \frac{1}{1-z}.$ We have  $ (A\ast B)(z) =A(z) 
\in \widetilde{TP_\infty}.$

\item We choose the sequence ${\mathbf B} =(b_k)_{k=0}^\infty  \in TP_\infty$ 
such that $b_k =k, B(z) = \frac{z}{(1-z)^2}.$ We have  $ (A\ast B)(z) = z A^\prime(z)  
\in \widetilde{TP_\infty}.$

\item We choose the sequence ${\mathbf B} =(b_k)_{k=0}^\infty  \in TP_\infty$ 
such that $b_k =k +1, B(z) = \frac{1}{(1-z)^2}.$ We have  $ (A\ast B)(z) =(zA(z))^\prime  
\in \widetilde{TP_\infty}.$

\item We choose the sequence ${\mathbf B} =(b_k)_{k=0}^\infty  \in TP_\infty$ 
such that $b_k = \frac{c^{k+1} - 1}{c -1}, c \in (0, 1) \cup (1, \infty), B(z) = 
\frac{1}{(1-z)} \cdot \frac{1}{(1- c  z)}.$ We have  $ (A\ast B)(z) =
\frac{1}{1 - c}(A(z) - c A(c z)) \in \widetilde{TP_\infty}.$

\item We choose the sequence ${\mathbf B} =(b_k)_{k=0}^\infty  \in TP_\infty$ 
such that $b_k = \frac{c^{k} - 1}{c -1}, c \in (0, 1) \cup (1, \infty), B(z) 
= \frac{z}{(1-z)} \cdot \frac{1}{(1-c  z)}.$ We have  $ (A\ast B)(z) =
\frac{1}{1 - c}(A(z) - A(c z))  \in \widetilde{TP_\infty}.$

\item We choose the sequence ${\mathbf B} =(b_k)_{k=0}^\infty  \in TP_\infty$ 
such that $ b_0 =1, b_k = d +1$  for $ k \geq 1, d >0, B(z) = 
\frac{1+ d z}{1-z}.$ We have  $ (A\ast B)(z) =(1+ d)(A(z)- \frac{d}{d +1}a_0)  
\in \widetilde{TP_\infty}.$

\item We choose the sequence ${\mathbf B} =(b_k)_{k=0}^\infty  \in TP_\infty$ 
such that $ b_k =0$ for $k=0, 1, \ldots, n,$ and $ b_k = 1$ for $ k \geq n+1,$
 $ B(z) = \frac{ z^{n+1}}{1-z}.$ We have  
$ (A\ast B)(z) =A(z)- \sum_{k=0}^n a_k z^k  \in \widetilde{TP_\infty}.$

\end{enumerate}
Lemma \ref{L1} is proved. $\Box$

By Lemma \ref{L1}(1), $A(z) \in \widetilde{TP_\infty},$  whence by theorem ASWE we have
\begin{equation}
A(z) =C z^q  e^{\gamma  z} \prod_{k=1}^\infty \frac{1+\alpha_k z}{1-\beta_k z},
\end{equation}
were $C >0, q\in \mathbb{N} \cup \{0\}, \gamma \geq 0,$ $\alpha_k \geq 0, \beta_k \geq 0,$
$\sum_{k=1}^\infty (\alpha_k + \beta_k) < \infty. $

By Lemma \ref{L1}(2), $A^\prime(z) \in \widetilde{TP_\infty}.$ By the assumption 
of Theorem~\ref{tt1}, the function  $A$  has at least one pole. Suppose  $A $  has at least
2 different positive poles. Since $A$ does not have positive zeros, then
$A^\prime$ has a positive zero (between poles), that is impossible. Thus,
$A$  has one (maybe, multiple) pole, and we have 
\begin{equation}
\label{e1}
A(z) =C z^q e^{\gamma  z}  \frac{\prod_{k=1}^\infty(1+\alpha_k z)}{(1-\beta z)^m}
=: \frac{F(z)}{(1-\beta z)^m} ,
\end{equation}
were $C >0, q\in \mathbb{N} \cup \{0\}, \gamma \geq 0,$ $\alpha_k \geq 0, \beta > 0,$
$\sum_{k=1}^\infty \alpha_k  < \infty. $

Since $F$ is an entire function with nonnegative Taylor coefficients,
we have $M(r, F) = \max_{|z| \leq r} |F(z)| =F(r).$ If $F$ is not a nonnegative constant,
then $\lim_{x\to +\infty} F(x) = + \infty.$ If $\lim_{x\to +\infty} 
\frac{F(x)}{(1-\beta x)^m} = + \infty,$ then $A^\prime$ has a positive zero
on $(\beta, +\infty),$ that is impossible. We conclude that
\begin{equation}
\label{e2}
A(z) =   \frac{C z^q \prod_{k=1}^n(1+\alpha_k z)}{(1-\beta z)^m}
=: \frac{P(z)}{(1-\beta z)^m} ,
\end{equation}
were $C >0, q\in \mathbb{N} \cup \{0\}, $ $\alpha_k \geq 0, \beta > 0,$
$n \in \mathbb{N}\cup \{0\}, q + n\leq m. $

Suppose that $m=2s, s\in \mathbb{N}.$  By Lemma \ref{L1}(5), $ (A\ast B)(z) =
\frac{1}{1 - c}(A(z) - A(c z))  \in \widetilde{TP_\infty}, c\neq 1.$ For $c>1$
we have $ (A\ast B)(z) = \frac{1}{c-1} \left (\frac{P(cz)}{(1-\beta c z)^{2s}} - 
\frac{P(z)}{(1-\beta z)^{2s}} \right).$ The function  $ (A\ast B)$ has 
two different poles: $\frac{1}{\beta c} < \frac{1}{\beta}. $ Since $P(x) >0$
for $x>0,$ and $c>1,$ we have $\lim_{x \to \frac{1}{\beta c} +0 } (A\ast B)(x) 
= +\infty, $ and $\lim_{x \to \frac{1}{\beta } -0 } (A\ast B)(x) 
= -\infty, $ so $ (A\ast B)$ has  a root on the interval $\left(\frac{1}{\beta c}, 
\frac{1}{\beta}\right).$ This contradicts to the fact that $A\ast B \in \widetilde{TP_\infty}.$
Thus, $m=2s+1, s\in \mathbb{N}\cup \{0\}.$

Suppose that $A(z) =\frac{P(z)}{(1-\beta z)^{2s+1}} $ with $\deg P =2s+1.$
Then $\lim_{x \to +\infty} A(x) = - L, L >0.$ By Lemma \ref{L1}(4), $ (A\ast B)(z) =
\frac{1}{c - 1}(cA(c z) - A(z))  \in \widetilde{TP_\infty}, c\neq 1.$ For $c>1$
we have  $ (A\ast B)(z) = \frac{1}{c-1} \left (\frac{c P(cz)}{(1-\beta c z)^{2s+1}} - 
\frac{P(z)}{(1-\beta z)^{2s+1}} \right).$ We observe that $\lim_{x \to \frac{1}{\beta } +0 } 
(A\ast B)(x) = + \infty, $  and $\lim_{x \to + \infty } 
(A\ast B)(x) = \frac{1}{c-1}(-cL +L) < 0. $  So, $ (A\ast B)$ has  a root on the interval 
$\left(\frac{1}{\beta }, 
+ \infty\right).$ This contradicts to the fact that $A\ast B \in \widetilde{TP_\infty}.$ Thus, 
$\deg P < 2s+1.$

We have proved that $A(z) =\frac{P(z)}{(1-\beta z)^{n}}, $ where $\deg P < n, n= 2s+1,
s\in \mathbb{N}\cup \{0\},$ and $P$ is a polynomial with nonnegative coefficients and all 
nonpositive roots. It remains to prove that $s=0.$ 

We observe that
\begin{equation}
\label{ef3}
A(z) =    \frac{B_0}{(1-\beta z)^n} -  \frac{B_1}{(1-\beta z)^{n-1}}
+ \ldots +(-1)^{n-1} \frac{B_{n-1}}{(1-\beta z)},
\end{equation}
where $P(z) =B_0 - B_1(1-\beta z) +  \ldots +(-1)^{n-1}(1-\beta z)^{n-1} =
B_0 + \beta B_1 (z- \frac{1}{\beta}) + \beta^2 B_2 (z- \frac{1}{\beta})^2 
+  \ldots + \beta^{n-1} B_{n-1} (z- \frac{1}{\beta})^{n-1},$ whence
$\beta^j B_j = \frac{P^{(j)}(\frac{1}{\beta})}{j!} >0$ for all $j=0, 1, \ldots, n-1.$

\begin{Lemma} \label{L2} Let $k \in \mathbb{N}, k \geq 2,$
 $A_{k, \beta}(z) = \frac{1}{(1-\beta z) ^k},$ and
 $F_{\gamma, \delta}(z)= \frac{e^{\gamma z}}{(1-\delta z)}, \gamma > 0, \delta >0.$ 
 Then 
 \begin{equation}
\label{e4}
(A_{k, \beta} \ast F_{\gamma, \delta})(z) =    
\frac{(k-1)!e^{\gamma \beta z}}{(1-\delta \beta z)^{k}} Q_{2k-2}(z),
\end{equation}
where $Q_{2k-2}(z)$ is a polynomial of degree at most $2k-2$ of the form
 \begin{equation}
\label{e5}
Q_{2k-2}(z) = \sum_{s=0}^{k-1} \sum_{t=k-s-1}^{k-1} 
\frac{ (\gamma \beta)^{k-1-s} (\delta \beta)^{t-k+s +1}}{(k-1-s)! (k-1-t)! t!} 
 z^t  (1-\delta \beta z)^{2k-t-s-2}.
\end{equation}
\end{Lemma}

{\it Proof of Lemma \ref{L2}.}  We have
$$ A_{k, \beta}(z) = \frac{1}{(k-1)! \beta^{k-1}} \left(\frac{1}{1-\beta z}\right)^{(k-1)} =
\frac{1}{(k-1)! \beta^{k-1}} \sum_{j=0}^\infty\beta^j (z^j)^{(k-1)}$$
$$= \frac{1}{(k-1)! \beta^{k-1}} \sum_{j=k-1}^\infty\beta^j j(j-1) \cdot \ldots 
\cdot (j-k+2) z^{j-k+1}     $$
$$  =   \frac{1}{(k-1)! \beta^{k-1}} \sum_{s=0}^\infty \beta^{s+k-1}(s+k-1)(s+k-2) 
\cdot \ldots \cdot (s+1) z^s$$
$$=  \frac{1}{(k-1)! } \sum_{s=0}^\infty \beta^{s}(s+k-1)(s+k-2) 
\cdot \ldots \cdot (s+1) z^s. $$
For a function $G(z) = \sum_{s=0}^\infty d_s z^s$ we obtain 
$$(A_{k, \beta} \ast G)(z) = \frac{1}{(k-1)! } \sum_{s=0}^\infty \beta^{s}(s+k-1)(s+k-2) 
\cdot \ldots \cdot (s+1) d_s  z^s  $$
$$= \frac{1}{(k-1)! } \left( z^{k-1} G(\beta z) \right)^{(k-1)}.   $$
Thus, for $G(z) = F_{\gamma, \delta}(z)$ we have
$$(A_{k, \beta} \ast F_{\gamma, \delta})(z) = \frac{1}{(k-1)! }
\left( e^{\gamma \beta z} z^{k-1}  (1-\delta \beta z)^{-1} \right)^{(k-1)}.   $$
Whence we get
$$(A_{k, \beta} \ast F_{\gamma, \delta})(z) =  \frac{1}{(k-1)! }\cdot $$
$$\sum_{s=0}^{k-1} \sum_{l=0}^s \frac{(k-1)!}{(k-1-s)! (s-l)! l!}
(e^{\gamma \beta z})^{(k-1-s)}(z^{k-1})^{(s-l)}((1-\delta \beta z)^{-1})^{(l)} $$
$$=  \sum_{s=0}^{k-1} \sum_{l=0}^s \left(\frac{1}{(k-1-s)! (s-l)! l!}
(\gamma \beta)^{k-1-s} e^{\gamma \beta z}\cdot \right. $$
$$\left. \frac{(k-1)!}{(k-s+l-1)!}z^{k-s+l-1} \cdot  l! (\delta \beta)^l 
\frac{1}{(1-\delta \beta z)^{l+1}}  \right) =$$
$$  \frac{(k-1)!e^{\gamma \beta z}}{(1-\delta \beta z)^{k}} 
\sum_{s=0}^{k-1} \sum_{l=0}^s \left(\frac{1}{(k-1-s)! (s-l)! (k-s+l-1)!} 
 \cdot \right. $$
$$\left. z^{k-s+l-1} (\gamma \beta)^{k-1-s} (\delta \beta)^l  (1-\delta \beta z)^{k-l-1}\right).$$
Changing the summation index $l$ in the second sum by  $k-s+l-1 =t,$ we obtain
$$(A_{k, \beta} \ast F_{\gamma, \delta})(z) =
 \frac{(k-1)!e^{\gamma \beta z}}{(1-\delta \beta z)^{k}} \cdot $$
$$ \left(\sum_{s=0}^{k-1} \sum_{t=k-s-1}^{k-1} \frac{(\gamma \beta)^{k-1-s} 
(\delta \beta)^{t-k+s+1}}{(k-1-s)! (k-1-t)! t!} \cdot 
 z^t   (1-\delta \beta z)^{2k-t-s-2}  \right)   $$
$$=:   \frac{(k-1)!e^{\gamma \beta z}}{(1-\delta \beta z)^{k}} Q_{2k-2}(z),   $$
where $Q_{2k-2}$ is a polynomial of degree at most $2k-2,$ since $\deg ( z^t 
(1-\delta \beta z)^{2k-t-s-2})
= 2k -s -2.$

Lemma \ref{L2} is proved. $\Box$

By (\ref{ef3}) we have
\begin{equation}
\label{e6}
A(z) =   B_0 A_{n, \beta}(z) - B_1 A_{n -1, \beta}(z) +   
\ldots +(-1)^{n-1} B_{n-1} A_{1, \beta}(z),
\end{equation}
where $B_0 >0, B_1>0, \ldots, B_{n-1} >0.$

So, using Lemma \ref{L2}, we get 
\begin{eqnarray}
\label{e7}
& \nonumber  (A \ast F_{\gamma, \delta})(z) =   B_0 (A_{n, \beta} \ast F_{\gamma, \delta})(z) - 
B_1 (A_{n -1, \beta}\ast F_{\gamma, \delta})(z)    
  \\ \nonumber &
+ \ldots +(-1)^{n-1} B_{n-1} (A_{1, \beta}\ast F_{\gamma, \delta})(z) =
\\ \nonumber &  B_0 \frac{(n-1)!e^{\gamma \beta z}}{(1-\delta \beta z)^{n}} Q_{2n-2}(z)
 -  B_1 \frac{(n-2)!e^{\gamma \beta z}}{(1-\delta \beta z)^{n-1}} Q_{2n-4}(z)
 \\  \nonumber & +  B_2 \frac{(n-3)!e^{\gamma \beta z}}{(1-\delta \beta z)^{n-2}} 
 Q_{2n-6}(z) - \ldots +(-1)^{n-1} B_{n-1} \frac{0!e^{\gamma \beta z}}{(1-\delta \beta z)} 
 Q_0(z) =
 \\ \nonumber & \frac{ e^{\gamma \beta z}}{(1-\delta \beta z)^n} \left(B_0(n-1)! 
 Q_{2n-2}(z) -  B_1(n-2)!  Q_{2n-4}(z)(1-\delta \beta z) + \right.
 \\ \nonumber & \left. B_2(n-3)!  Q_{2n-6}(z)(1-\delta \beta z)^2  - \ldots +(-1)^{n-1} 
B_{n-1} 0! Q_0(z) (1-\delta \beta z)^{n-1} \right)
\\  & =: \frac{ e^{\gamma \beta z}}{(1-\delta \beta z)^n} H_{2n-2}(z),
\end{eqnarray}
where the degree of a polynomial $H_{2n-2}$ is at most $2n-2.$

Since $F_{\gamma, \delta}(z)= \frac{e^{\gamma z}}{(1-\delta z)} \in \widetilde{TP_\infty}$
for all $ \gamma > 0, \delta >0,$ by our assumption we conclude that
$ (A \ast F_{\gamma, \delta}) \in \widetilde{TP_\infty},$ whence the polynomial
$H_{2n-2}$ has all nonnegative coefficients and all nonpositive roots.
We denote by
$$ H_{2n-2}(z) =:\sum_{s=0}^{2n-2} h_s z^s, h_s \geq 0. $$
For $n >1 $ we have $2n-3 > 0,$ and we want to evaluate $h_{2n-3.}$
Note that  $z^{2n-3}$ can be found only in the terms $B_0(n-1)! 
 Q_{2n-2}(z)$ and $ B_1(n-2)!  Q_{2n-4}(z)(1-\delta \beta z)$  of formula  (\ref{e7}). By (\ref{e5}),
$$Q_{2n-2}(z) = \sum_{s=0}^{n-1} \sum_{t=n-s-1}^{n-1} 
\frac{ (\gamma \beta)^{n-1-s} (\delta \beta)^{t-n+s+1}}{(n-1-s)! (n-1-t)! t!} 
 z^t  (1-\delta \beta z)^{2n-t-s-2}.$$
We observe that $\deg ( z^t (1-\delta \beta z)^{2n-t-s-2})
= 2n - s -2 < 2n-3$ for $s>1,$ so we will search for the term $z^{2n-3}$
in the summands with $s=0$ and $s=1.$ We have 
$$Q_{2n-2}(z) = \frac{ (\gamma \beta)^{n-1} (\delta \beta)^{0}}{((n-1)!)^2 } 
 z^{n-1}  (1-\delta \beta z)^{n-1}  +  $$
 $$ \frac{ (\gamma \beta)^{n-2} (\delta \beta)^{0}}{((n-2)!)^2 } 
 z^{n-2}  (1-\delta \beta z)^{n-1}   +
 \frac{ (\gamma \beta)^{n-2} (\delta \beta)^{1}}{(n-2)! (n-1)! } 
 z^{n-1}  (1-\delta \beta z)^{n-2}$$
  $$+  \sum_{s=2}^{n-1} \sum_{t=n-s-1}^{n-1} 
\frac{ (\gamma \beta)^{n-1-s} (\delta \beta)^{t-n+s+1}}
{(n-1-s)! (n-1-t)! t!}  z^t  (1-\delta \beta z)^{2n-t-s-2}.$$
Thus, gathering the terms with $z^{2n-3}$ in the first 3 summands
of the above formula, we obtain the term with $z^{2n-3}$ in $B_0(n-1)! 
 Q_{2n-2}(z):$ 
$$ B_0(n-1)! \left(\frac{ (\gamma \beta)^{n-1} (\delta \beta)^{0}}{((n-1)!)^2 }
(-1)^{n-2}  (n-1)(\delta \beta)^{n-2} +\right.$$
$$\left.\frac{ (\gamma \beta)^{n-2} (\delta \beta)^{0}}{((n-2)!)^2 } 
(-1)^{n-1} (\delta \beta)^{n-1}   +   
\frac{ (\gamma \beta)^{n-2} (\delta \beta)^{1}}{(n-2)! (n-1)! } 
(-1)^{n-2} (\delta \beta)^{n-2}\right) =$$
$$B_0  \frac{(-1)^{n-2} (\gamma \beta)^{n-2} (\delta \beta)^{n-2}}{(n-2)! }
( \gamma \beta - \delta \beta (n-1) + \delta \beta) $$
$$= B_0\frac{(-1)^{n-2} (\gamma \beta)^{n-2} (\delta \beta)^{n-2}}{(n-2)! }
( \gamma \beta - \delta \beta (n-2)).$$

By (\ref{e5}),
$$  Q_{2n-4}(z)(1-\delta \beta z) =$$
$$ (1-\delta \beta z) \cdot
 \sum_{s=0}^{n-2} \sum_{t=n-s-2}^{n-2} 
\frac{ (\gamma \beta)^{n-2-s} (\delta \beta)^{t-n+s +2}}{(n-2-s)! (n-2-t)! t!} 
 z^t  (1-\delta \beta z)^{2n-t-s-4}.$$
 We observe that $\deg ( z^t (1-\delta \beta z)^{2n-t-s-4})
= 2n - s -4 < 2n-4$ for $s\geq 1,$ so we will search for the term $z^{2n-4}$
in the summands with $s=0$ and multiply it by $(-\delta \beta z).$  Thus, the
 term with $z^{2n-3}$ in $ (-B_1(n-2)!  Q_{2n-4}(z)(1-\delta \beta z))$ equals
$$- B_1(n-2)! \frac{(-1)^{n-1}(\gamma \beta)^{n-2}(\delta \beta)^{n-1} }{((n-2)!)^2}  
= B_1 \frac{(-1)^{n}(\gamma \beta)^{n-2}(\delta \beta)^{n-1} }{(n-2)!}.$$
Finally, we get 
$$ h_{2n-3} = B_0\frac{(-1)^{n-2} (\gamma \beta)^{n-2} (\delta \beta)^{n-2}}{(n-2)! }
( \gamma \beta - \delta \beta (n-2)) + 
B_1 \frac{(-1)^{n}(\gamma \beta)^{n-2}(\delta \beta)^{n-1} }{(n-2)!}$$
$$ = \frac{(-1)^{n-2} (\gamma \beta)^{n-2} (\delta \beta)^{n-2}}{(n-2)! }(B_0
\gamma \beta -B_0(n-2) \delta \beta + B_1\delta \beta).$$

Since $ n= 2s+1, s\in \mathbb{N}\cup \{0\},$ and $\beta >0, \gamma >0, \delta >0,$
we obtain
$$\mathrm{sign} h_{2n-3} = - \mathrm{sign} (B_0
\gamma \beta -B_0(n-2) \delta \beta - B_1\delta \beta)   = -1 $$
for $\gamma >0$ being large enough and $\delta >0$ being small enough.
We get a contradiction. Thus, $n=1$ and $A(z) = \frac{C}{1- \beta z}, C > 0, \beta > 0.$

Theorem \ref{tt1} is proved. 

\section{Proof of Theorem \ref{tt2}}

Let us prove the necessity.  Let ${\mathbf A} =(a_k)_{k=0}^\infty$ 
be a nonnegative sequence, such that $a_k>0$ for $0\leq k \leq 3,$  and 
$a_k=0$ for $k\geq 4.$ Suppose that the operator $\Lambda_{\mathbf A}$
preserves the set of the $TP_\infty-$sequences. Since the sequence
${\bf B}_1 =(1, 1, 1, 1, \ldots) \in TP_\infty,$ we have 
 $\Lambda_{\mathbf A} ({\bf B}_1) = (a_0, a_1, a_2, a_3, 0, 0, 0, \dots)
 \in TP_\infty,$ whence the polynomial $\sum_{k=0}^3 a_k x^k$
 has only real (and nonpositive) zeros. Further, since the sequence
${\bf B}_2 =(0, 1, 1, 1, 1, \ldots) \in TP_\infty,$ we have 
 $\Lambda_{\mathbf A} ({\bf B}_2) = (0, a_1, a_2, a_3, 0, 0, 0, \dots)
 \in TP_\infty,$ whence the polynomial $\sum_{k=1}^3 a_k x^k$
 has only real (and nonpositive) zeros. The necessity is proved.
 
Let us prove the sufficiency. Obviously, $(c_k)_{k=0}^{\infty} \in TP_{\infty}$ if 
and only if $(C \lambda ^k c_k)_{k=0}^{\infty} \in TP_{\infty}$ for $C>0, \ \lambda >0.$ 
Thus, without loss of generality we can assume that $a_0=a_1=1.$  Then we can rewrite 
our sequence $\mathbf{A}$ in the form $\left(1,1,\frac{1}{a}, \frac{1}{a^2 b},0,0,0,\ldots \right),$ 
where $a=\frac{a_1^2}{a_0 a_2}=\frac{1}{a_2}, \ b=\frac{a_2^2}{a_1 a_3}=\frac{a_2^2}{a_3}.$ 

So, by assumption both the polynomial

\begin{equation}
\label{e1}
P(x, a, b) = 1 + x + \frac{x^2}{a} + \frac{x^3}{a^2 b}, \  a>0, b> 0.
\end{equation}
and the polynomial
\begin{equation}
\label{e2}
T(x, a, b) =  x + \frac{x^2}{a} + \frac{x^3}{a^2 b}.
\end{equation}
have only real non-positive zeros.

Note that
\begin{equation}
\label{e3}
P(x, a, b) = 1 +T(x, a, b).
\end{equation}
Denote by
\begin{equation}
\label{e4}
F(y, a, b) = 1 +ay(1 + y +\frac{y^2}{b})
\end{equation}
and by
\begin{equation}
\label{e5}
t(y,  b) = y + y^2 +\frac{y^3}{b}= y(1 + y +\frac{y^2}{b}).
\end{equation}
We have
\begin{equation}
\label{e6}
F(y, a, b)= 1 + a t(y,  b).
\end{equation}

\begin{Statement} \label{s1}  Both $P$ and $T$
have only real zeros if and only if both $F$ and $t$
have only real zeros.
\end{Statement}

{\it Proof.} Statement \ref{s1} follows from the two identities below.
\begin{equation}
\label{e7}
F(y, a, b)= P(ay, a, b)
\end{equation}
and
\begin{equation}
\label{e8}
t(y,  b) = \frac{1}{a} T(ay, a, b). 
\end{equation}
$\Box$

The following fact is obvious.
\begin{Statement} \label{s2}  The polynomials $T(x,a,b)$ and 
$t(y,  b)$ have only real zeros if and only if $b \geq 4.$
\end{Statement}

From now on we will assume that 
\begin{equation}
\label{g1}
b \geq 4.
\end{equation}
Consider the derivative of the polynomial $F(y,a,b).$
\begin{equation}
\label{g2}
F'(y,a,b)=a\left(1+2y+\frac{3}{b} y^2\right)
\end{equation}
Denote by
\begin{equation}
\label{g3}
\alpha_1(b):=\frac{b}{3}\left(-1+\sqrt{1-\frac{3}{b}} \right)
\end{equation}
and
\begin{equation}
\label{g4}
\alpha_2(b):=\frac{b}{3}\left(-1-\sqrt{1-\frac{3}{b}} \right)
\end{equation}
the roots of $F'(y,a,b).$
It follows from (\ref{g1}) that
\begin{equation}
\label{g5}
\frac{\alpha_1(b)}{b}=\frac{1}{3}\left(-1+
\sqrt{1-\frac{3}{b}} \right) \geq-\frac{1}{6}.
\end{equation}
Since all roots of $F(y,a,b)$ are real, $\alpha_2(b) < \alpha_1(b),$ 
and $\lim_{y\to+\infty}F(y,a,b)=+\infty,$ it is clear that
\begin{equation}
\label{g6}
F\left(\alpha_1(b),a,b\right)\leq 0,
\end{equation}
and that
\begin{equation}
\label{g7}
F\left(\alpha_2(b),a,b\right)\geq 0.
\end{equation}

The inequality (\ref{g7}) can be improved in the following way.
\begin{Statement}  \label{s3}
\begin{equation}
\label{g8}
F\left(\alpha_2(b),a,b\right)\geq 1.
\end{equation}
\end{Statement}

{\it Proof.} It follows from (\ref{e5}) that
\begin{equation}
\label{g9}
t'(y,b)=\frac{1}{a}F'(y,a,b).
\end{equation}
Therefore, $\alpha_1(b)$ from (\ref{g3}) and $\alpha_2(b)$ from (\ref{g4}) are roots of 
$t'(y,b)$ too. By our assumption, $t(y,b)$ has only real roots. Since
$\lim_{y\to-\infty}t(y,b)=-\infty,$ we have
\begin{equation}
\label{g10}
t\left(\alpha_2(b),b\right)\geq 0
\end{equation}
By virtue of (\ref{e6}), the following is true
$$F\left(\alpha_2(b),a,b\right)=1+a t\left(\alpha_2(b),b\right)\geq 1.$$
Statement \ref{s3} is proved. $\Box$

In notations introduced in (\ref{e1}) and (\ref{e2}) the sufficiency in Theorem~\ref{tt2} 
could be equivalently reformulated in the following form.
\begin{Statement} \label{t1}  Assume that both $P(x, a, b)$
and $T(x, a, b)$ have only real non-positive
zeros, and $G(x) =\sum_{k=0}^\infty c_k x^k\in \widetilde{TP}.$ 
Then $P\ast G \in \widetilde{TP}.$
\end{Statement}
First, consider the case when $c_0=0.$ If additionally $c_1=0,$ the 
 Statement \ref{t1} is obvious.

Let $c_1\neq 0.$ Since $G(x)\in \widetilde{TP}$, and therefore, $G(x)\in \widetilde{TP}_2,$ 
we conclude that
\begin{equation}
\label{g11}
c_2^2-c_1c_3\geq 0.
\end{equation}
Since $T(x,a,b)$ has only real zeros we have $b\geq 4$ (see (\ref{g1})). Thus,
\begin{equation}
\label{g12}
b c_2^2-4c_1c_3\geq 0.
\end{equation}
The last inequality means that the polynomial
$$P\ast G (x,a,b)=c_1x+\frac{c_2}{a}x^2+\frac{c_3}{a^2b}x^3$$
has only real roots, that is by ASWE theorem $P\ast G \in \widetilde{TP}.$ So, in this case Statement \ref{t1} is true.

From now on we will assume that $c_0 \neq0.$
For $G(x) =\sum_{k=0}^\infty c_k x^k$
we denote by 
\begin{equation}
\label{e20}
p_k =\frac{c_{k-1}}{c_k},\    k \in\mathbb{N},
\end{equation}
and by
\begin{equation}
\label{e21}
q_{k+1} =\frac{p_{k+1}}{p_k},\    k \in\mathbb{N}.
\end{equation}
Then we have
\begin{eqnarray}
\label{e22} &
G(x) = c_0 \left(1+ \frac{x}{p_1}+\frac{x^2}{p_1p_2} + \ldots +
\frac{x^k}{p_1 p_2\cdot \ldots \cdot p_k } + \ldots \right) =
\\  \nonumber &
 c_0 \left(1+ \frac{x}{p_1} +\left(\frac{x}{p_1}\right)^2 \cdot\frac{1}{q_2}
+ \left(\frac{x}{p_1}\right)^3 \cdot\frac{1}{q_2^2 q_3}+ \ldots + 
\left(\frac{x}{p_1}\right)^k \cdot\frac{1}{q_2^k q_3^{k-1}\cdot \ldots 
\cdot q_{k-1}^2 q_k}+\ldots \right).
\end{eqnarray}
Denote by
\begin{equation}
\label{e24}
g(y) = 1+ y +\frac{y^2}{q_2} + \frac{y^3}{q_2^2 q_3}
+ \ldots + \frac{y^k}{q_2^k q_3^{k-1}\cdot \ldots 
\cdot q_{k-1}^2 q_k} + \ldots,
\end{equation}
so that by (\ref{e22})
\begin{equation}
\label{e25}
g(y) = \frac{1}{c_0} G( p_1 y)\in \widetilde{TP}.
\end{equation}

\begin{Statement} \label{s7}  If $g(y)= 1 + y +
\sum_{k=2}^\infty \frac{y^k}{q_2^{k-1}q_3^{k-2}\cdot \ldots \cdot q_k}
\in \widetilde{TP},$  then 
\begin{equation}
\label{e46} 
q_2 \geq 1, \ \ q_3 \geq 1, \ \ q_2^2 q_3 - 2 q_2 q_3 + 1\geq 0.
\end{equation}
\end{Statement}

{\it Proof.} By the definition, the statement $g(y) \in \widetilde{TP}$
means that all minors of the matrix
\begin {equation}
\label{e46a} 
 \left[
  \begin{array}{cccccc}
   1 & 1 & \frac{1}{q_2} & \frac{1}{q_2^2 q_3} & \frac{1}{q_2^3 q_3^2 q_4} &\ldots \\
   0   &1 & 1 & \frac{1}{q_2}& \frac{1}{q_2^2 q_3} &\ldots \\
   0   &  0  & 1 & 1& \frac{1}{q_2} &\ldots \\
      \vdots&\vdots&\vdots&\vdots&\vdots&\vdots
  \end{array}
 \right]
\end {equation}
are nonnegative. The first statement of (\ref{e46}) is equivalent to the
fact that
$$ \left|
  \begin{array}{cc}
   1 & \frac{1}{q_2}  \\
   1   &1  \\
     \end{array}
 \right|  \geq 0.$$
The second statement of (\ref{e46}) follows from the
fact that
$$ \left|
  \begin{array}{cc}
   \frac{1}{q_2} & \frac{1}{q_2^2 q_3}  \\
   1   & \frac{1}{q_2}  \\
     \end{array}
 \right|  \geq 0.$$
The third statement of (\ref{e46}) can be easily obtained from the inequality below
$$ \left|
  \begin{array}{ccc}
1&   \frac{1}{q_2} & \frac{1}{q_2^2 q_3}  \\
 1&  1   & \frac{1}{q_2}  \\
 0 & 1 & \frac{1}{q_2}\\
     \end{array}
 \right|  \geq 0.$$
Statement \ref{s7} is proved.  $\Box$

By (\ref{e1}) we have
\begin{eqnarray}
\label{e23}
&
(P\ast G) (x, a, b)=c_0 \left(1 +\frac{x}{p_1} +
 \left(\frac{x}{p_1}\right)^2  \cdot\frac{1}{q_2 a}
\right. \\ \nonumber &
\left.  + \left(\frac{x}{p_1}\right)^3 \cdot
\frac{1}{(q_2 a)^2 q_3 b}\right).
\end{eqnarray}
Denote by
\begin{equation}
\label{e26}
F_q (y, a, b) =1+q_2 a y +q_2 a y^2+\frac{q_2 a}{q_3 b} y^3  =1 + q_2 a y \left(1 + y + 
\frac{y^2}{q_3 b}\right).
\end{equation}
By (\ref{e23}) we have
\begin{equation}
\label{e27}
F_q (y, a, b) = \frac{1}{c_0} (P\ast G) (y p_1q_2 a, a, b).
\end{equation}
Now we can equivalently reformulate Statement \ref{t1}
in the following way.
\begin{Statement} \label{t2}  Assume that both $F(y, a, b)$ and
$t(y, b)$ from (\ref{e4}) and (\ref{e5}) have only real non-positive zeros,
and $g(y)$ from (\ref{e24})  belongs to $\widetilde{TP}.$ Then $F_q (y, a, b) \in \widetilde{TP}. $
\end{Statement}

{\it Proof.} We will show that 
\begin {equation}
\label{ef47}   F_q\left(\frac{\alpha_1(b)}{q_2}, a, b\right) \leq 0, 
\end{equation}
and
\begin {equation}
\label{ef48}   F_q(\alpha_2(b), a, b) \geq 0.
\end{equation}
Since $\lim_{x\to +\infty}F_q(x,a,b)=+\infty$ and $\lim_{x\to -\infty}F_q(x,a,b)=-\infty,$ 
Statement \ref{t2} will follow from (\ref{ef47}) and (\ref{ef48}).

Let us prove (\ref{ef47}).  It follows from (\ref{e26}) that
\begin {equation}
\label{ef49}   
F_q\left(\frac{\alpha_1(b)}{q_2}, a, b\right)=1+a\alpha_1(b)+\frac{a}{q_2}\left(\alpha_1(b)\right)^2+
\frac{a}{bq_2^2q_3}\left(\alpha_1(b)\right)^3,
\end{equation}
and from (\ref{e4}) that
\begin {eqnarray}
\label{ef50}   
& F_q(\alpha_1(b), a, b) -F_q\left(\frac{\alpha_1(b)}{q_2}, a, b\right)=
\\ \nonumber &  a\left(\alpha_1(b)\right)^2\left(1-\frac{1}{q_2}\right)+
\frac{a}{b}\left(\alpha_1(b)\right)^3\left(1-\frac{1}{q_2^2q_3}\right).
\end{eqnarray}
By (\ref{g5}) we have
\begin {equation}
\label{ef51}   
 F_q(\alpha_1(b), a, b) -F_q\left(\frac{\alpha_1(b)}{q_2}, a, b\right)\geq 
\end{equation}
$$a\left(\alpha_1(b)\right)^2\left(\left(1-\frac{1}{q_2}\right)-\frac{1}{6}\left(1-\frac{1}{q_2^2q_3}\right)\right)=$$
$$\frac{a\left(\alpha_1(b)\right)^2}{6q_2^2q_3}\left(5q_2^2q_3-6q_2q_3+1\right).$$
Note that by the first statement of (\ref{e46}) we obtain
\begin{equation}
\label{ef52}   
 \left(5q_2^2q_3-6q_2q_3+1\right)-\left(q_2^2q_3-2q_2q_3+1\right)=
 4q_2q_3(q_2-1) \geq 0.
\end{equation}
It follows from (\ref{ef51}) and the third statement of (\ref{e46}) that
\begin {eqnarray}
\label{ef53}  
& 
 F_q(\alpha_1(b), a, b) -F_q\left(\frac{\alpha_1(b)}{q_2}, a, b\right)\geq
 \\  \nonumber & \frac{a\left(\alpha_1(b)\right)^2}{6q_2^2q_3}
 \left(q_2^2q_3-2q_2q_3+1\right)\geq 0.
\end{eqnarray}
Thus, by virtue of (\ref{g6}) we conclude that 
$$F_q\left(\frac{\alpha_1(b)}{q_2}, a, b\right) \leq F_q(\alpha_1(b), a, b)\leq 0.$$
So, (\ref{ef47}) is proved.

Let us prove (\ref{ef48}).
It follows from (\ref{e26}) and the second statement of (\ref{e46}) that 
$$  F_q(\alpha_2(b), a, b) = 1+ q_2 a \alpha_2(b)
\left(1+\alpha_2(b) + \frac{\left(\alpha_2(b)\right)^2}{q_3 b}\right) >$$
$$ 1+ q_2 a \alpha_2(b) \left(1+\alpha_2(b) + 
\frac{\left(\alpha_2(b)\right)^2}{ b}\right). $$
By (\ref{e4}) and (\ref{g8}) we have
\begin {equation}
\label{e48}   F_q(\alpha_2(b), a, b) > 1 + q_2(F(\alpha_2(b), a, b)
-1) \geq 1.
\end{equation}
So, (\ref{ef48}) and thereby Statement \ref{t2} is proved.$\Box$

Theorem \ref{tt2}  is proved. 

{\bf Aknowledgement.} \thanks{ The authors are indebted to Professor
Alan Sokal for posing the problem and for the inspiring discussions during 
the AIM workshop ``Theory and applications of total positivity'',
July 24-July 28, 2023 (see \cite{conf} for more details). 
We are also very grateful to the American Institute of Mathematics and to the 
organizers of this remarkable workshop Shaun Fallat, Dominique Guillot, and 
Apoorva Khare for the possibility to participate in the workshop, to meet the 
colleges and to discuss many interesting and important questions on total positivity.


\begin{thebibliography}{20}

\bibitem{aissen}
M. Aissen, A. Edrei, I.J. Schoenberg, A. Whitney, On the
Generating Functions of Totally Positive Sequences, \textsl{ J.
Anal. Math.} \textbf{2} (1952), 93--109.

\bibitem{ando}
T. Ando, Totally positive matrices,  \textsl{ Linear Alg. \& Its Appl.} , 90 
(1987),  165--219

\bibitem{cc1}
T. Craven and G. Csordas,  Complex zero decreasing sequences,
\textsl{Methods Appl. Anal.}, \textbf{2} (1995), 420--441.

\bibitem{Dya1}
A. Dyachenko, Total Nonnegativity of Infinite Hurwitz Matrices of Entire and
Meromorphic Functions,  \textsl{Complex Analysis and Operator Theory}, 
 \textbf{8}, (2014), 1097--1127.

\bibitem{Dya2}
A. Dyachenko, Hurwitz matrices of doubly infinite series, \textsl{ Linear Algebra 
and its Applications}, \textbf{530}, (2017), 266--287.

\bibitem{DyaSok} 
A. Dyachenko and A. Sokal, Total-positivity characterization of the
Laguerre-P\'olya class $LP^+$. Pages 661--662 in: K. Driver, O. Holtz
and A. Sokal (eds.), Oberwolfach Rep. 19 (2022), no. 1, 657--681, DOI:
10.4171/OWR/2022/13

\bibitem{fek}
M. Fekete, G. P\'olya, \"Uber ein Problem von Laguerre, \textsl{
Rend. Circ. Mat. Palermo} \textbf{34} (1912), 89--120.

\bibitem{HW}
I. I. Hirschman and D.V.Widder, {\em  The Convolution Transform},
Princeton University Press, Princeton, New Jersey, 1955.

\bibitem{hut}
J. I. Hutchinson,  On a remarkable class of entire functions,
\textsl{Trans. Amer. Math. Soc.} \textbf{25} (1923), 325--332.

\bibitem{tp}
S. Karlin,  \textsl{Total Positivity,} Vol. I, Stanford University
Press, California 1968.

\bibitem{klv}
O.M.Katkova, T.Lobova, A.M.Vishnyakova, On power series
having sections with only real zeros,  \textsl{ Computation
Methods and Functional Theory}, \textbf{3}, No 2, (2003), 425--441.

\bibitem{klv1}
O.M.Katkova, T.Lobova, A.M.Vishnyakova, On entire functions having 
Taylor sections with only real zeros,  \textsl{ Journal of Mathematical Physics, 
Analysis, Geometry}, \textbf{11}, No. 4, (2004), 449--469. 

\bibitem{kosshap}
V.P.Kostov, B.Shapiro, 
Hardy-Petrovitch-Hutchinson's problem and 
partial theta function,  \textsl{Duke Math. J.}, \textbf{162}, No. 5 (2013), 
825--861. 

\bibitem{lev}
B. Ja. Levin, \textsl{Distribution of Zeros of Entire Functions,}
Transl. Math. Mono., 5, Amer. Math. Soc., Providence, RI, 1964;
revised ed. 1980.

\bibitem{ngthv1}
Thu Hien Nguyen, Anna Vishnyakova, On the entire functions from the Laguerre--P\'olya 
class having the decreasing second quotients of Taylor coefficients, \textsl{Journal of 
Mathematical Analysis and Applications}, \textbf{465}, No. 1 (2018), 348 -- 359.

\bibitem{HishAn}
Thu Hien Nguyen, Anna Vishnyakova, Hutchinson's intervals and entire 
functions from the Laguerre-P\'olya class,  arXiv:2212.05692v1 [math.CV].

\bibitem{O}
N.~Obreschkov , Verteilung und Berechnung der Nullstellen reeller
Polynome, VEB Deutscher Verlag der Wissenschaften, Berlin, 1963.

\bibitem{iv}
I. V. Ostrovskii,  On Zero Distribution of Sections and Tails of
Power Series, \textsl{Israel Math. Conference Proceedings,}
\textbf{15} (2001), 297--310.

\bibitem{ivzh}
I.V. Ostrovskii, N.A. Zheltukhina, Parametric representation of a class of 
multiply positive sequences, \textsl{ Complex Variables, Theory and Application,}
\textbf{ 37}, 1-4, (1998), 457--469.

\bibitem{pin}
A. Pinkus, Spectral properties of totally positive kernels and matrices, 
M. Gasca (ed.) C.A. Micchelli (ed.) , Total Positivity and its Applications, 
Kluwer Acad. Publ. (1996), 477--511.

\bibitem{pol}
G. P\'olya, Collected Papers, Vol. II Location of Zeros, (R.P.Boas
ed.) MIT Press, Cambridge, MA, 1974.

\bibitem{polsch}
G. P\'olya and J.Schur, \"Uber zwei Arten von Faktorenfolgen in der Theorie
der algebraischen Gleichungen, J. Reine Andrew. Math., 144 (1914),
pp. 89-113.

\bibitem{War} S.~O.~Warnaar, Partial theta functions, 
\url{ https://www.researchgate.net/publication/327791878_Partial_theta_functions}.


\bibitem{conf}
The AIM workshop "Theory and applications of total positivity", 
 \url{https://aimath.org/pastworkshops/totalpos.html}


\end{thebibliography}
\end{document}